\documentclass[a4paper,12pt]{article}
\usepackage[top=2.5cm,bottom=2.5cm,left=2.5cm,right=2.5cm]{geometry}
\usepackage{amsfonts}
\usepackage{latexsym,amsmath,amssymb,amsfonts,epsfig,psfrag,url,graphics,ifpdf,multicol}
\usepackage{mathtools} 
\usepackage{color,colordvi,amscd,amsthm} 
\usepackage{tikz}
\usepackage{verbatim}
\usepackage[numbers,sort&compress]{natbib}
\usepackage{indentfirst,graphics,epsfig}
\usepackage{graphicx}
\usepackage{graphics}
\usepackage{graphicx,psfrag}
\usepackage{graphics}
\usepackage{mathrsfs}
\usepackage{tabularx}
\usepackage{booktabs}
\usepackage{floatrow}
\usepackage{multirow}
\usepackage{graphicx}
 \usepackage{enumitem}
\usepackage{caption}
\usepackage[ruled]{algorithm2e}
\newfloatcommand{capbtabbox}{table}[][\FBwidth]

\newtheorem{thm}{Theorem}

\theoremstyle{plain}
\theoremstyle{plain}

 \theoremstyle{definition}

\newtheorem{prop}{Proposition}[section]
\theoremstyle{remark}
\def\pf{\noindent{\bf \textbf{Proof}.\ }}
\def\qed{{\hfill\rule{4pt}{7pt}}}

\allowdisplaybreaks

\numberwithin{subcase}{case}
\numberwithin{subca}{ca}
\numberwithin{subsubcase}{subcase}

\usepackage{indentfirst,subfig}
\usepackage{booktabs} 
\begin{document}
	\captionsetup[figure]{labelfont={bf},name={Fig.},labelsep=period}

	\begin{center} {\large The stability of independence polynomials of complete bipartite graphs}
	\end{center}
	\pagestyle{empty}
	
	\begin{center}
		{
			{\small Guo Chen$^a$, Bo Ning$^b$, Jianhua Tu$^{a,}$\footnote{Corresponding author.\\\indent \ \  E-mail: chenguo030603@163.com (G. Chen), bo.ning@nankai.edu.cn (B. Ning), tu@btbu.edu.cn (J. Tu)}}\\[2mm]
			
			{\small $^a$School of Mathematics and Statistics, Beijing Technology and Business University, \\
				\hspace*{1pt} Beijing, 100048, PR China \\
				$^b$College of Cryptology and Cyber Science, Nankai University, Tianjin, 300350, PR China}\\[2mm]}
		
	\end{center}
	
	\baselineskip=0.24in
		\begin{center}
		\begin{abstract}
	     The independence polynomial of a graph is termed {\it stable} if all its roots are located in the left half-plane $\{z \in \mathbb{C} : \mathrm{Re}(z) \leq 0\}$, and the graph itself is also referred to as stable. Brown and Cameron (Electron. J. Combin. 25(1) (2018) \#P1.46) proved that the complete bipartite graph $K_{1,n}$ is stable and posed the question: \textbf{Are all complete bipartite graphs stable?}
	     We answer this question by establishing the following results:
	  	\begin{itemize}
		
		\item The complete bipartite graphs $K_{2,n}$ and $K_{3,n}$ are stable.  
		
		\item For any integer $k\geq0$, there exists an integer $N(k)\in \mathbb{N}$ such that $K_{m,m+k}$ is stable for all $m>N(k)$.

         \item For any rational $\ell> 1$, there exists an integer $N(\ell) \in \mathbb{N}$ such that whenever $m >N(\ell)$ and $\ell \cdot m$ is an integer, $K_{m, \ell \cdot m}$ is \textbf{not} stable.
			
		\end{itemize}
			
		\end{abstract}
	\end{center}

	\section{Introduction}\label{sec1}
	
	Let $G$ be a simple graph. An independent set of $G$ is a subset of vertices of $G$ such that no two vertices in this subset are adjacent. The independence number $\alpha(G)$ of $G$ is the cardinality of a maximum independent set in $G$. The independence polynomial of $G$, denoted by $i(G,x)$, is the generating polynomial that counts the number of independent sets of each size in $G$, and is defined by  
	\[i(G,x) = \sum_{k=0}^{\alpha(G)} i_k(G) x^k,\]
	where $i_k(G)$ is the number of independent sets of size $k$ in $G$. The roots of $i(G,x)$ are called the independence roots of $G$. For decades, the study of independence polynomials and independence roots has been a vibrant research area, with several significant open problems remaining unresolved. See, for example, \cite{Brown2005, Brown2018, Chudnovsky2007, Csikvari2013, Gutman1983, Levit2005, Wang2011, Zhang2022, Zhu2018} for related work on independence polynomials and their roots.
	
	The celebrated result due to Chudnovsky and Seymour \cite{Chudnovsky2007} establishes that all independence roots of a claw-free graph are real. Since the coefficients of the independence polynomial are positive integers, it follows that all independence roots of a claw-free graph lie on the negative real axis. Furthermore, it is known \cite{Brown2000} that for any graph, the independence root with the smallest modulus is always located on the negative real axis. On the other hand, Brown and Nowakowski \cite{Brown2005} have demonstrated that almost all graphs have a non-real independence root. It has also been observed \cite{Brown2004} that for small graphs, all independence roots lie in the left half-plane $\{z\in\mathbb{C}: \mathrm{Re}(z) \leq 0\}$. Given these findings, the following questions naturally emerge: 
		\vskip 3mm
	{\it  
	How prevalent are graphs that have all their independence roots in the left half-plane?}
	\vskip 3mm
	
	This question is interesting as it aim to further our understanding of the properties and distribution of independence roots of graphs. The answers could potentially enhance our knowledge of the relationship between the structure of a graph and the location of its independence roots.
	
	A polynomial whose roots all lie in the left half-plane is termed {\it Hurwitz quasi-stable} (or simply {\it stable}). Such polynomials hold significant importance in numerous applied fields \cite{Choe2004}. A graph is also referred to as stable if its independence polynomial is stable. Brown and Cameron \cite{Brown2018} showed that every graph $G$ with $\alpha(G) \leq 3$ is stable. However, they also proved that for any graph 
	$G$ with $\alpha(G) \geq 4$, there exists a graph $H$ with $\alpha(H)=\alpha(G)$ such that $G$ is an induced subgraph of $H$ and $H$ is not stable.
	Additionally, they established the following result:
		
	\begin{thm}\cite{Brown2018}\label{thm1}
		The complete bipartite graph $K_{1,n}$ is stable.
	\end{thm}
		
	Calculations show that small complete bipartite graphs are stable, so Brown and Cameron \cite{Brown2018} posed the question: \textbf{Are all complete bipartite graphs stable?} This paper shall answer this question. For the complete bipartite graph $K_{m,n}$ (without loss of generality, assume $m\leq n$), we first analyze cases with small $m$ and establish:
		
	\begin{thm}\label{thm2}
		The complete bipartite graphs $K_{2,n}$ and $K_{3,n}$ are stable.
	\end{thm}
	
	Next, we demonstrate the stability of $K_{m,n}$ with nearly balanced bipartitions where the size difference is bounded by a constant, provided that $m$ is sufficiently large.
	
	\begin{thm}\label{thm3}
	Let $k\in \mathbb{N}$ be an integer. 
	\begin{itemize}
	\item For any integer $k\leq 6$, $K_{m,m+k}$ is stable for all $m\geq 1$. 
	
	\item For any integer $k\geq 7$, there exists an integer $N(k)\in\mathbb{N}$ such that $K_{m,m+k}$ is stable for all $m> N(k)$. 
	\end{itemize}
	\end{thm}
	
	Finally, we prove that for the complete bipartite graph $K_{m,n}$, if the sizes of its two parts differ by a fixed factor $\ell>1$ (i.e., the parts are unbalanced in size), then $K_{m,n}$ is \textbf{not} stable when $m$ is sufficiently large.
	
		\begin{thm}\label{thm4}
     For any rational $\ell> 1$, there exists an integer $N(\ell) \in \mathbb{N}$ such that whenever $m >N(\ell)$ and $\ell\cdot m$ is an integer, $K_{m, \ell\cdot m}$ is \textbf{not} stable.
		\end{thm}

	The proofs of Theorems \ref{thm2}, \ref{thm3}, and \ref{thm4} will be presented in Sections \ref{sec2}, \ref{sec3}, and \ref{sec4}, respectively.

		
	\section{Proof of Theorem \ref{thm2}}\label{sec2}
	
	The independence polynomial of the complete bipartite graph $K_{m,n}$ is 
	\[i(K_{m,n},x)=(1+x)^n+(1+x)^m-1.\]
	Let $y=1+x$, and define
	\[\Phi(K_{m,n},y)=y^n+y^m-1.\]
	Thus, $K_{m,n}$ is stable if and only if every root $z$ of $\Phi(K_{m,n}, y)$ satisfies $\mathrm{Re}(z) \leq 1$.
	
	When $m=n$, the polynomial becomes $\Phi(K_{m,m}, y)=2 y^m-1$. Its roots are 
	\[y_k = \sqrt[m]{\frac{1}{2}} \cdot e^{i \cdot \frac{2k\pi}{m}}, \quad k = 0, 1, 2, \ldots, m-1,\]  
	which lie on the circle centered at the origin with radius $\sqrt[m]{\frac{1}{2}}$. Since this radius is less than 1, all roots satisfy $\mathrm{Re}(z) \leq 1$. Consequently, $K_{m,m}$ is stable.  
	
	We introduce the well-known Rouch\'{e}'s Theorem (see \cite{Fisher1990}), which serves as a key tool in proving Theorems \ref{thm2} and \ref{thm3}.
	
	\begin{thm}[Rouch\'{e}'s Theorem]\label{thmR}  
		Let $f(z)$ and $g(z)$ be analytic functions inside and on a simple piecewise smooth closed curve $\gamma$. If  
		\[|f(z) + g(z)| < |f(z)| \quad \text{for all } z \in \gamma, \]  
		then $f(z)$ and $g(z)$ have the same number of zeros (counting multiplicities) in the interior of $\gamma$.  
	\end{thm}

	We present the following proposition that will be used in proving that $K_{2,n}$ and $K_{3,n}$ are stable.		

    \begin{prop}\label{prop2.1}
	The complete bipartite graph $K_{m,m+1}$ is stable.
    \end{prop}

    \pf We will use Rouch\'{e}'s Theorem to prove that if $z$ is a root of $\Phi(K_{m,m+1},y)$, then $\mathrm{Re}(z)<1$.
	
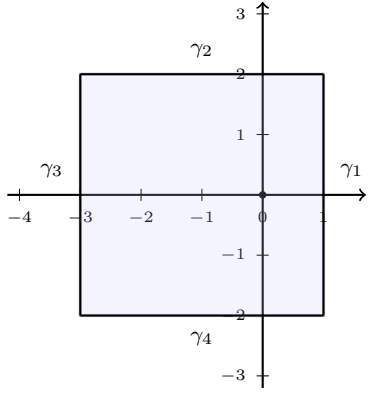
\begin{figure}[h]
				\centering
				\begin{tikzpicture}[scale=0.8]
					
					\draw[->, thick, black] (-4.2,0) -- (1.7,0); 
					\draw[->, thick, black] (0,-3.2) -- (0,3.2); 
					
					\foreach \x in {-4,-3,-2,-1,0,1}
					\draw (\x,0.1) -- (\x,-0.1) node[below, font=\tiny]{$\x$};
					
					\foreach \y in {-3,-2,-1,1,2,3}
					\draw (0.1,\y) -- (-0.1,\y) node[left, font=\tiny]{$\y$};
					\draw (0.1,0) -- (-0.1,0);
					\filldraw[black] (0,0) circle (1.5pt);
					
					\filldraw[fill=blue!20, draw=black, very thick, opacity=0.2] 
					(-3,-2) rectangle (1,2);
					
					\draw[thick, black] (1,-2) -- (1,2) node[pos=0.6, right=2pt, font=\scriptsize] {$\gamma_1$};
					\draw[thick, black] (-3,-2) -- (-3,2) node[pos=0.6, left=2pt, font=\scriptsize] {$\gamma_3$}; 
					\draw[thick, black] (-3,2) -- (1,2) node[pos=0.5, above=2pt, font=\scriptsize] {$\gamma_2$};
					\draw[thick, black] (-3,-2) -- (1,-2) node[pos=0.5, below=2pt, font=\scriptsize] {$\gamma_4$};
				\end{tikzpicture}
				\caption{The curve $\gamma$ and the closed region it bounds.}
				\label{fig1}
			\end{figure}

	Let $f(z)=-z^{m+1}-z^{m}$, $g(z)=\Phi(K_{m,m+1},z)=z^{m+1}+z^{m}-1$ and set
	\begin{itemize}
	\item  $\gamma_1=\left\lbrace z\mid \mathrm{Re}(z)=1  \ \text{and} \ -2\leq \mathrm{Im}(z) \leq 2\right\rbrace;$			
			
	\item $\gamma_2=\left\lbrace z\mid -3 \leq \mathrm{Re}(z) \leq 1 \  \text{and} \ \mathrm{Im}(z) =2\right\rbrace;$
			
	\item $\gamma_3=\left\lbrace z\mid \mathrm{Re}(z)=-3  \  \text{and} \ -2\leq \mathrm{Im}(z) \leq 2\right\rbrace;$
			
	\item  $\gamma_4=\left\lbrace z\mid -3 \leq \mathrm{Re}(z) \leq 1 \  \text{and} \ \mathrm{Im}(z) =-2\right\rbrace.$		
	\end{itemize}
	Let $\gamma$ be the curve composed of four line segments $\gamma_1, \gamma_2, \gamma_3, \gamma_4$, as illustrated in Figure \ref{fig1}.
			
		We now show that $|f(z)+g(z)|<|f(z)|$ for all $z\in \gamma$. Note that the functions $f(z)$ and $g(z)$ are analytic on the entire complex plane, with
		$|f(z)+g(z)|=1$ and $|f(z)|=|z|^m\cdot|z+1|$. Thus, we need to show that 
		\[|z|^m\cdot |z+1|>1 \quad \text{for all\ } z\in \gamma.\]
		
		It is observed that 
		\begin{itemize}
			\item[(1)] for any $z\in \gamma$, $|z| \geq 1$, with equality if and only if $z=1$,
			
			\item[(2)] for any $z\in \gamma$, $|z+1| \geq 2$, with equality if and only if $z\in\{1,-3,-1+2i,-1-2i\}$.
		\end{itemize}
		Therefore, \[|z|^m\cdot |z+1|>1 \quad \text{for all\ } z\in \gamma.\] 
		By Rouch\'{e}'s Theorem, $f(z)=-z^{m+1}-z^{m}$ and $g(z)=z^{m+1}+z^{m}-1$ have the same number of zeros inside $\gamma$ counting multiplicities. We know that $f(z)$ has a simple root at $z=-1$ and a root at $z=0$ of multiplicity $m$. Thus, $g(z)=\Phi(K_{m,m+1},z)$ has all $m+1$ of its roots in the interior of $\gamma$, which implies that if $z$ is a root of $\Phi(K_{m,m+1},y)$, then $\mathrm{Re}(z)<1$. 
		
		The proof of Proposition \ref{prop2.1} is complete. \qed
		
		
       
      
       \begin{prop}\label{prop2.2}
      	The complete bipartite graph $K_{2,n}$ is stable.
      \end{prop}
		
	   \pf We establish the stability of the complete bipartite graph $K_{2,n}$ by demonstrating that if $z$ is a root of $\Phi(K_{2,n},x)$, then $\mathrm{Re}(z)<1$. By Proposition \ref{prop2.1}, $K_{2,3}$ is stable. Thus, we assume that $n\geq 4$.
	   
	   Let $f(z)=-z^n$ and $g(z)=\Phi(K_{2,n},z)=z^n+z^2-1$. Let $\gamma$ be the curve introduced in the proof of Proposition \ref{prop2.1} and shown in Figure \ref{fig1}.
	   We now show that $|f(z)+g(z)|<|f(z)|$ for all $z\in \gamma$. Note that the functions $f(z)$ and $g(z)$ are analytic on the entire complex plane, with
	   $|f(z)+g(z)|=|z^2-1|=|z+1|\cdot|z-1|$ and $|f(z)|=|z|^n$. Thus, we need to show that 
	   \[|z|^n> |z+1|\cdot|z-1|\quad \text{for all\ } z\in \gamma.\]
	   
        We distinguish the following three cases.
        
        \textbf{Case 1:} $z\in \gamma_1$.
			 
		Let $z=1+ki$ where $-2 \leq k \leq 2$. Now, $|z|=\sqrt{k^2+1}$, $|z-1|=|k|$, and $|z+1|=\sqrt{k^2+4}$. Since 
		\[\left( \sqrt{k^2+1}\right) ^n \geq \left( \sqrt{k^2+1}\right) ^4=\sqrt{k^8+4k^6+6k^4+4k^2+1}>\sqrt{k^4+4k^2}=|k|\cdot\sqrt{k^2+4},\]
		it follows that $|z|^n> |z+1|\cdot|z-1|$ for all $z\in \gamma_1$.
			
		\textbf{Case 2:} $z\in \gamma_2\cup \gamma_4$.
		
		Let $z=k\pm2i$ where $-3 \leq k \leq 1$. Now $|z|=\sqrt{4+k^2}$, $|z+1|=\sqrt{4+(k+1)^2}$, and $|z-1|=\sqrt{4+(k-1)^2}$. It follows that 
		\begin{align*}
			\lvert z + 1 \rvert \cdot \lvert z - 1 \rvert 
			&= \sqrt{4 + (k + 1)^2} \cdot \sqrt{4 + (k - 1)^2} \\
			&\leq \sqrt{4 + 2^2} \cdot \sqrt{4 + 4^2}= \sqrt{160} \\
			&< \left( \sqrt{k^2 + 4} \right)^n = \lvert z \rvert^n.
		\end{align*}
		Hence $|z|^n> |z+1|\cdot|z-1|$ for all $z\in \gamma_2\cup \gamma_4$.
			
		\textbf{Case 3:} $z\in \gamma_3$.
			
	Let $z=-3+ki$ where $-2 \leq k \leq 2$. Now $|z|=\sqrt{k^2+9}$, $|z+1|=\sqrt{k^2+4}$, and $|z-1|=\sqrt{k^2+16}$. It follows that
	\[|z+1|\cdot|z-1|=\sqrt{k^4+20k^2+64}<\sqrt{(k^2+9)^4}=\left( \sqrt{k^2+9}\right) ^4 \leq |z|^n.\]
	Hence $|z|^n> |z+1|\cdot|z-1|$ for all $z\in \gamma_3$.
			
    Now, we have shown that for all $z \in \gamma$, $|f(z)+g(z)|<|f(z)|$. Thus, by Rouch\'{e}'s Theorem, we obtain that $f(z)=-z^n$ and $g(z)=z^n+z^2-1$ have the same number of zeros inside $\gamma$. Since $f(z)$ has a root at $z=0$ of multiplicity $n$, $g(z)=\Phi(K_{2,n},z)$ has all $n$ of its roots in the interior of $\gamma$, which implies that if $z$ is a root of $\Phi(K_{2,n},y)$, then $\mathrm{Re}(z)<1$. 
    
    The proof of Proposition \ref{prop2.2} is complete. \qed

\begin{figure}[h]
    	\centering
    	\begin{tikzpicture}[scale=0.8]
    		\draw[->, thick, black] (-4.2,0) -- (0.7,0); 
    		\draw[->, thick, black] (0,-3.2) -- (0,3.2); 
    		\foreach \x in {-4,-3,-2,-1}
    		\draw (\x,0.1) -- (\x,-0.1) node[below, font=\tiny]{$\x$};
    		\foreach \y in {-3,-2,-1,0,1,2,3}
    		\draw (0.1,\y) -- (-0.1,\y) node[left, font=\tiny]{$\y$};
    		\filldraw[black] (0,0) circle (1.5pt);
    		\filldraw[fill=blue!20, draw=black, very thick, opacity=0.2] 
    		(-3,-2) rectangle (0,2);
    		
    			\draw[thick, black] (0,-2) -- (0,2) node[pos=0.6, right=2pt, font=\scriptsize] {$\beta_1$};
    		\draw[thick, black] (-3,-2) -- (-3,2) node[pos=0.6, left=2pt, font=\scriptsize] {$\beta_3$}; 
    		\draw[thick, black] (-3,2) -- (0,2) node[pos=0.5, above=2pt, font=\scriptsize] {$\beta_2$};
    		\draw[thick, black] (-3,-2) -- (0,-2) node[pos=0.5, below=2pt, font=\scriptsize] {$\beta_4$};
    	\end{tikzpicture}
    	\caption{The curve $\beta$ and the closed region it bounds.}
    	\label{fig2}
    \end{figure}
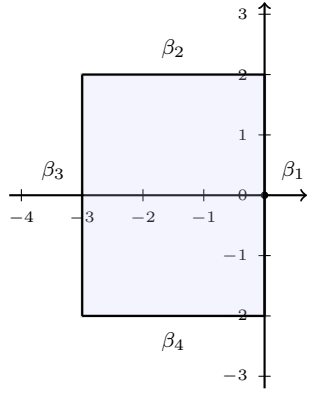

    Next, we will prove that $K_{3,n}$ is stable. The proof will utilize results on the distribution of roots of $\Phi(K_{1,n}, y)$.     
    Set 
    	\begin{itemize}
    	\item  $\beta_1=\left\lbrace z\mid \mathrm{Re}(z)=0  \ \text{and} \ -2\leq \mathrm{Im}(z) \leq 2\right\rbrace;$			
    	
    	\item $\beta_2=	\left\lbrace z\mid -3 \leq \mathrm{Re}(z) \leq 0 \  \text{and} \ \mathrm{Im}(z) =2\right\rbrace;$
    	
    	\item $\beta_3=	\left\lbrace z\mid \mathrm{Re}(z)=-3  \  \text{and} \ -2\leq \mathrm{Im}(z) \leq 2\right\rbrace;$
    	
    	\item  $\beta_4=\left\lbrace z\mid -3 \leq \mathrm{Re}(z) \leq 0 \  \text{and} \ \mathrm{Im}(z) =-2\right\rbrace.$		
    \end{itemize}
  Let $\beta$ be the curve composed of four line segments $\beta_1, \beta_2, \beta_3, \beta_4$, as illustrated in Figure \ref{fig2}. It was demonstrated by Brown and Cameron \cite{Brown2018} that all roots of $i(K_{1,n},x)$ lie in the interior of $\beta$. Consequently, $\Phi(K_{1,n},y) = y^n + y - 1$ has all its roots contained in the interior of $\gamma$, which is depicted in Figure \ref{fig1}.  
    
     \begin{prop}\label{prop2.3}
    	The complete bipartite graph $K_{3,n}$ is stable.
    \end{prop}
    
    \pf We establish the stability of the complete bipartite graph $K_{3,n}$ by demonstrating that if $z$ is a root of $\Phi(K_{3,n},x)$, then $\mathrm{Re}(z)<1$. By Proposition \ref{prop2.1}, $K_{3,4}$ is stable. Thus, we assume that $n\geq 5$.
    
    Let $f(z)=-\Phi(K_{1,n},z)=-z^n-z+1$ and $g(z)=\Phi(K_{3,n},z)=z^n+z^3-1$. Let $\gamma$ be the curve introduced in the proof of Proposition \ref{prop2.1} and shown in Figure \ref{fig1}.
    We now show that $|f(z)+g(z)|<|f(z)|$ for all $z\in \gamma$. Note that the functions $f(z)$ and $g(z)$ are analytic on the entire complex plane, with
    $|f(z)+g(z)|=|z^3-z|$ and $|f(z)|=|z^n+z-1|$.

    \textbf{Case 1:} $z\in \gamma_2\cup \gamma_3\cup \gamma_4$.
    
    In this case, $|z|\geq2$. By the triangle inequality, we have 
    \[|f(z)+g(z)|=|z^3-z|\leq |z|^3+|z|,\]
    and
    \[|f(z)|=|z^n+z-1| \geq |z|^n-|z|-1 \geq |z|^5-|z|-1.\]
    Thus,
    \[|f(z)|-|f(z)+g(z)|\geq |z|^5-|z|^3-2|z|-1.\]
    Let $p(x)=x^5-x^3-2x-1$. The function $p(x)$ is strictly increasing on $[2,+\infty)$, which implies that for any $x\in [2,+\infty)$, $p(x)\geq p(2)=19$. Hence, $|f(z)+g(z)|<|f(z)|$ for all $z\in \gamma_2\cup \gamma_3\cup \gamma_4$.
    
    \textbf{Case 2:} $z\in \gamma_1$.
    
	Let $z=1+ki$ where $-2 \leq k \leq 2$. Now we have 
    \[|f(z)+g(z)|=|z^3-z|=|z|\cdot|z-1|\cdot|z+1|=\sqrt{k^2+1}\cdot |k| \cdot \sqrt{k^2+4},\]
    and 
    \[|f(z)|=|z^n+z-1|=|z^n+ki|\geq |z|^n-|k|\geq \left(\sqrt{k^2+1} \right)^5-|k|.\]
    
   Let $q(x) = x^5 + 5x^4 + 9x^3 + 3x^2 - 5x + 1$ for $x \in [0, 4]$. The function $q(x)$ is convex on this interval and attains its global minimum value of $0.05297$ at $x \approx 0.3066$. Thus,
    \begin{align*}
    	&\left(\sqrt{k^2+1}\cdot |k| \cdot \sqrt{k^2+4}+|k|\right)^2\\
    	=&\left(|k|\cdot(1+\sqrt{k^2+1}\sqrt{k^2+4})\right)^2\\
    	=& k^2(1+k^4+5k^2+4+2\sqrt{k^2+1}\sqrt{k^2+4})\\
    	\leq &k^2(1+k^4+5k^2+4+(k^2+1)+(k^2+4))\\
    	=&k^6+7k^4+10k^2\\
    	<&(1+k^2)^5,
    \end{align*}
    which implies that $|f(z)+g(z)|<|f(z)|$ for all $z\in \gamma_1$.

  We have shown that $|f(z) + g(z)| < |f(z)|$ for all $z \in \gamma$. By Rouch\'{e}'s Theorem, combined with the fact that $f(z) = -\Phi(K_{1,n}, z)$ has all its roots inside $\gamma$, it follows that $g(z) = \Phi(K_{3,n}, z)$ also has all $n$ of its roots in the interior of $\gamma$. Consequently, if $z$ is a root of $\Phi(K_{3,n}, y)$, then $\mathrm{Re}(z) < 1$.
    
   The proof of Proposition \ref{prop2.3} is complete. \qed

   \noindent{\it Proof of Theorem \ref{thm2}}. The theorem follows immediately by combining Propositions \ref{prop2.2} and \ref{prop2.3}. \qed

	
	\section{Proof of Theorem \ref{thm3}}\label{sec3}
	
     \noindent{\it Proof of Theorem \ref{thm3}}. Since in Section \ref{sec2} we have shown that Theorem \ref{thm3} holds when $k=0$ and $k=1$, we assume $k \geq 2$. 
     
     Let $k \geq 2$ be a given integer, define $f(z)=-z^{m+k}-z^m$ and $g(z)=\Phi(K_{m,m+k},z)=z^{m+k} + z^m - 1$. Let $\gamma$ be the curve introduced in the proof of Proposition \ref{prop2.1} and shown in Figure \ref{fig1}. Note that the functions $f(z)$ and $g(z)$ are analytic on the entire complex plane, with
     $|f(z)+g(z)|=1$ and $|f(z)|=|z|^m \cdot |z^k+1|$. 
     
     If $z\in \gamma_2\cup \gamma_3\cup \gamma_4$, then $|z|\geq 2$. Thus, 
     \[|f(z)|=|z|^m \cdot |z^k+1|\geq |z|^m\cdot (|z|^k-1)>1.\]
     Hence $|f(z)+g(z)|<|f(z)|$ for all $z\in \gamma_2\cup \gamma_3\cup\gamma_4$.
		
    We consider the case when $z \in \gamma_1$. Let  
    \[z = 1 + \tan\theta \cdot i = \sec\theta (\cos\theta + i\sin\theta),\]  
    where $\theta \in [-\arctan 2, \arctan 2]$. Note that $\sec\theta \in [1, \sqrt{5}]$. Thus,  
      \begin{align*}
    	|f(z)|^2 &= \left(|z|^m \cdot |z^k + 1|\right)^2 \\
    	&= \left(|z|^m \cdot \left|\sec^k\theta (\cos k\theta + i\sin k\theta) + 1\right|\right)^2 \\
    	&= \sec^{2m}\theta \cdot \left[ \left(\sec^k\theta \cdot \cos k\theta + 1\right)^2 + \left(\sec^k\theta \cdot \sin k\theta\right)^2 \right] \\
    	&= \sec^{2m}\theta \cdot \left( \sec^{2k}\theta + 1 + 2\sec^k\theta \cos k\theta \right).
    \end{align*}

    Using the identities $\cos(0\theta) = 1$, $\cos(1\theta) = \cos\theta$, and  
    \[\cos(k\theta) = 2\cos\theta \cos((k-1)\theta) - \cos((k-2)\theta),\]  
    it follows that $\cos k\theta$ is a polynomial in $\cos\theta$ of degree at most $k$. Consequently, $2\sec^k\theta \cos k\theta$ is a polynomial in $\sec\theta$. Let  
    \[P_{k}(\sec\theta) = \sec^{2k}\theta + 1 + 2\sec^k\theta \cos k\theta,\]  
    and let $t=\sec\theta$. Then $P_k(t)$ is a polynomial in $t$ of degree $2k$ defined on $[1, \sqrt{5}]$. Let $c_k$ denote the minimum value of $P_k(t) $ over this interval, i.e.,  
    \[c_k = \min_{t \in [1, \sqrt{5}]} P_k(t).\]  
    Direct computation yields 
    \begin{align*}
    	&P_2(t) = t^4 - 2t^2 + 5, \quad c_2 = 4, \\
    	&P_3(t) = t^6 - 6t^2 + 9, \quad c_3 \approx 3.3431, \\
    	&P_4(t) = t^8 + 2t^4 - 16t^2 + 17, \quad c_4 \approx 2.3555, \\
    	&P_5(t) = t^{10} + 10t^4 - 40t^2 + 33, \quad c_5 \approx 1.6123, \\
    	&P_6(t) = t^{12} - 2t^6 + 36t^4 - 96t^2 + 65, \quad c_6 \approx 1.1249.
    \end{align*}
        
    Therefore, for any positive integer $k \leq 6$,  
    \begin{equation}\label{eq1}
    	|f(z)|^2= \sec^{2m}\theta \cdot P_k(\sec\theta) \geq \sec^{2m}\theta \cdot c_k > 1,
    \end{equation}
     which implies $|f(z) + g(z)| < |f(z)|$ for all $z \in \gamma_1$.
    
    Thus, for $k \leq 6$, $|f(z) + g(z)| < |f(z)|$ for all $z \in \gamma$. By Rouch\'{e}’s Theorem, combined with the fact that $f(z)$ has all its roots inside $\gamma$, it follows that $g(z) = \Phi(K_{m,m+k}, z)$ also has all $m+k$ roots in the interior of $\gamma$. Consequently, if $z$ is a root of $\Phi(K_{m,m+k}, y)$, then $\mathrm{Re}(z) < 1$. This proves that $K_{m,m+k}$ is stable for any integer $k\leq 6$ and all $m\geq 1$.

   
   Calculations reveal that $c_k>1$ does not necessarily hold when $k \geq 7$. For example, when $k = 7$,  
   \[P_7(t) = t^{14} - 14t^6 + 112t^4 - 224t^2 + 129, \quad \text{with } c_7 \approx 0.8105.\]  
   Thus, for $k \geq 7$, we cannot guarantee that Inequality (\ref{eq1}) holds for all $m\geq 1$. However, we can prove that for any fixed $k \geq 7$, there exists an integer $N(k)\in\mathbb{N}$ such that the inequality $|f(z)|^2>1$ holds for all $z \in \gamma_1$ and all $m > N(k)$. To establish this, we first analyze the properties of $P_k(t)$.
   
   If $\theta = 0$, then $z = 1$ and $P_k(\sec\theta) = |z^k + 1|^2 = 4$, implying $P_k(1) = 4$. For $\theta \neq 0$, we have $|z| = \sec\theta > 1$, which implies that $ z^k + 1 \neq 0$ and $P_k(\sec\theta) = |z^k + 1|^2 > 0$. Therefore, $P_k(t) > 0$ for all $t \in [1, \sqrt{5}]$, which ensures $c_k > 0$.
   
   Fix $k \geq 7$. We may assume $c_k \leq 1$, as otherwise Inequality (\ref{eq1}) would hold. There exists $\delta(k) > 0$ such that $P_k(t) \geq 2$ for all $t \in [1, 1 + \delta(k)]$. Define  
   \[N(k) =\left\lceil\frac{\ln\frac{1}{c_k}}{2\ln(1 + \delta(k))}\right\rceil.\]  
   
   For $\sec\theta \in [1, 1 + \delta(k)]$,  
   \[|f(z)|^2 = \sec^{2m}\theta \cdot P_k(\sec\theta) \geq \sec^{2m}\theta \cdot 2 > 1.\]  
   When $\sec\theta \in [1 + \delta(k), \sqrt{5}]$ and $m > N(k)$,  
   \[|f(z)|^2 = \sec^{2m}\theta \cdot P_k(\sec\theta) > (1 + \delta(k))^{\ln\frac{1}{c_k}/\ln(1 + \delta(k))} \cdot c_k = (1 + \delta(k))^{\log_{1+\delta(k)}\frac{1}{c_k}} \cdot c_k = 1.\]  
   Thus, for $m >N(k) $, we have $|f(z) + g(z)| < |f(z)|$ for all $z \in \gamma_1$.
   
   Consequently, for any fixed $k \geq 7$, when $m > N(k)$, Rouch\'{e}’s Theorem implies that $g(z) = \Phi(K_{m,m+k}, z)$ has all $m + k$ roots in the interior of $\gamma$, as $f(z)$ already has all its roots inside $\gamma$. Therefore, if $z$ is a root of $\Phi(K_{m,m+k}, y)$, then $\mathrm{Re}(z) < 1$. This proves that for any fixed $k \geq 7$, $K_{m,m+k}$ is stable for all $m > N(k)$.
   
   The proof of Theorem \ref{thm3} is complete. \qed

	
	\section{Proof of Theorem \ref{thm4}}\label{sec4}

	\noindent{\it Proof of Theorem \ref{thm4}.} Let $\ell=\frac{p}{q}>1$ be a fixed rational number  where $p$ and $q$ are coprime positive integers. If $\ell\cdot m$ is an integer, then $q\mid m$. Let $m=q\cdot r$. Thus,  
	\[
	\Phi(K_{m,\ell\cdot m}, y) = y^{pr} + y^{qr} - 1.
	\]
	Let $z = y^r$ and define 
	\[
	g(z) := z^p + z^q - 1.
	\]
	We note that $g(z)$ must have a root $t$ with $|t| > 1$. This is because of Vieta's formulas, which show that for all roots $z_1, \dots, z_p$ of $g(z)$,
	\[
	\prod_{i=1}^p z_i = (-1)^{p+1}.
	\]
	Hence 
	\begin{equation}\label{eq2}
		\prod_{i=1}^p |z_i| = 1.
	\end{equation}
	Observe that $g(0) = -1 < 0$ and $g(1) = 1 > 0$. By the Intermediate Value Theorem, there is a real root $\xi \in (0,1)$, so $|\xi| < 1$. From Equation (\ref{eq2}), since the product of the moduli of all roots is $1$ and one root has modulus less than $1$, there must be at least one root $t$ with $|t| > 1$.
	
	Let $t$ be the unique root of $g(z)$ satisfying:  
	\begin{itemize}
	\item[1.] $|t|$ is maximum among all roots,  
	
	\item[2.] Among roots with maximum modulus, $\mathrm{Re}(t)$ is largest,  
	
	\item[3.] If multiple roots still satisfy (1) and (2), choose the one with $\mathrm{Im}(t)>0$.  
	\end{itemize}
	The preceding analysis demonstrates $|t|>1$. 
	Write 
	$t=s\cdot e^{i\theta},$
	where $s,\theta\in \mathbb{R}$, $s>1$, and $\theta \in \left[ 0,2\pi \right) $. Thus, $\Phi(K_{m,\ell \cdot m},y)$ has $r$ roots of the form:
		\[y_k = \sqrt[r]{s} \cdot e^{i\frac{\theta + 2k\pi}{r}}, \quad k = 0,1,2,\cdots,r-1.\]
	Note that both $s$ and $\theta$ are functions of $\ell$.

	\textbf{Case 1: $\theta\in[0,\pi]$.}
	
    In this case, we will show that there exists an integer $U(s,\theta)$ such that $\mathrm{Re}(y_0)>1$ for all $r>U(s,\theta)$.
	
	If $\theta=0$, then 
	\[\mathrm{Re}(y_0)=\sqrt[r]{s}>1,\]
	and let $U(s,\theta)=0$.
	If $\theta \in (0,\pi]$, then 
	 \[\mathrm{Re}(y_0)=\sqrt[r]{s}\cdot \cos(\theta/r).\]
	
	Consider the function $f(t)=s^t\cdot \cos(t\theta)$ for $t\in (0,1/2)$. We have
		\[f'(t)=s^t\cdot\left[\ln s\cdot \cos(t\theta)-\theta\cdot \sin(t\theta)\right].\]
	Let 
   \[\delta(s,\theta)=\min\left\lbrace \frac{1}{2},\frac{\arctan(\frac{\ln s}{\theta})}{\theta}\right\rbrace.\] 
   Since $h(t)=\ln s\cdot \cos(t\theta)-\theta\cdot \sin(t\theta)$ is strictly decreasing on $(0,1/2)$, we have
   \[f'(t)>0\quad \text{for all}\ t\in(0,\delta(s,\theta)).\]
   Thus, $f(t)$ is strictly increasing on $(0,\delta(s,\theta))$.
   
   As
   \[\lim_{t\to 0^+}f(t)=\lim_{t\to 0^+}s^t\cdot \cos(t\theta)=1,\] 
   it follows that for any $t\in (0,\delta(s,\theta))$, $f(t)>1$. Let
   \[U(s,\theta)=\left\lceil \frac{1}{\delta(s,\theta)}\right\rceil.\]
   Then, 
   \[\mathrm{Re}(y_0)=\sqrt[r]{s}\cdot \cos(\theta/r)>1 \quad \text{for all}\ r>U(s,\theta).\]

   	\textbf{Case 2: $\theta\in(\pi,2\pi)$.}
   
   In this case,
   \[y_{r-1}=\sqrt[r]{s} \cdot e^{i\frac{\theta + 2(r-1)\pi}{r}}=\sqrt[r]{s} \cdot e^{i\frac{\theta-2\pi}{r}}\]
   which implies that 
   \[\mathrm{Re}(y_{r-1})=\sqrt[r]{s}\cdot \cos(\frac{2\pi-\theta}{r}).\]
   Let $\theta'=2\pi-\theta$, then $\theta'\in (0,\pi)$. Thus, 
   \[\mathrm{Re}(y_{r-1})=\sqrt[r]{s}\cdot \cos(\theta'/r)>1 \quad \text{for all}\ r>U(s,\theta').\]

 Let 
 \[
 N(\ell):= \begin{cases}
 	q\cdot U(s,\theta), & 0 \leq \theta \leq \pi \\
 	q\cdot U(s,2\pi-\theta), & \pi < \theta < 2\pi.
 \end{cases}
 \]
 Thus, whenever $m > N(\ell)$ and $\ell\cdot m$ is an integer, there exists a root $z$ of $\Phi(K_{m,\ell\cdot m}, y)$ such that $\mathrm{Re}(z)>1$, which implies that $K_{m,\ell\cdot m}$ is not stable.

   The proof of Theorem \ref{thm4} is complete. \qed


	\section*{Use of AI tools declaration}
	The authors declare that they have not used Artificial Intelligence (AI) tools in the creation of this article.
	
	\section*{Data Availability}
	Data sharing is not applicable to this paper as no datasets were generated or analyzed during
	the current study.
	
	\section*{Acknowledgments}
	This work was supported by Beijing Natural Science Foundation (No. 1232005). 
	
	\section*{Conflict of interest}
	The authors declare that they have no conflicts of interest.

\end{document}